\documentclass[twoside,reqno,A4]{amsart}

\allowdisplaybreaks
\usepackage{latexsym,rotate,eucal,cite,url}
\usepackage{amsmath,amsthm,amssymb,amsxtra} 

\newcommand{\quab}{\hspace*{2.0mm}}
\newcommand{\puab}{\hspace*{-1.5mm}}

\newcommand{\pav}[1]{\lfloor{#1}\rfloor}

\newcommand{\pavv}[1]{\left\lfloor{#1}\right\rfloor}

\newcommand{\ang}[1]{\langle{#1}\rangle}

\newcommand{\mb}[1]{\mathbb{#1}}

\newcommand{\ff}[3]{\left(#1;#2\right)_{#3}}
\newcommand{\hh}[3]{\left[#1;#2\right]_{#3}}
\newcommand{\hyp}[3]{\left[\puab\ba{#1}#2\\#3\ea\puab\right]}
\newcommand{\hyq}[4]{\left[\puab\ba{#1}#3\\#4\ea{\puab\Big|#2}\right]}

\newcommand{\binm}{\binom}
\newcommand{\binq}[2]{{#1\brack#2}}

\newcommand{\nnm}{\nonumber}
\newcommand{\be}{\begin{equation}}
\newcommand{\ee}{\end{equation}}
\newcommand{\ba}{\begin{array}}
\newcommand{\ea}{\end{array}}
\newcommand{\bmn}{\begin{eqnarray}}
\newcommand{\emn}{\end{eqnarray}}
\newcommand{\bnm}{\begin{eqnarray*}}
\newcommand{\enm}{\end{eqnarray*}}
\newcommand{\bln}{\begin{subequations}}
\newcommand{\eln}{\end{subequations}}

\newcommand{\pq}[1]{\begin{equation}#1\end{equation}}
\newcommand{\pmq}[1]{\begin{align}#1
            \end{align}}     
\newcommand{\pnq}[1]{\begin{align*}#1
            \end{align*}}    
\newcommand{\pmp}[2]{\begin{alignat}{#1}#2
            \end{alignat}}   
\newcommand{\mult}[2]{\begin{array}{#1}#2\end{array}}%

\newcommand{\centro}[1]
           {\begin{center}#1\end{center}}

\newcommand{\alp}{\alpha}
\newcommand{\bet}{\beta}
\newcommand{\gam}{\gamma}

\newcommand{\lam}{\lambda}
\newcommand{\vph}{\varphi}

\newcommand{\Lam}{\Lambda}

\newcommand{\Gam}{\Gamma}

\newtheorem{thm}{Theorem}

\newtheorem{corl}[thm]{Corollary}
\newtheorem{prop}[thm]{Proposition}

\newtheorem{entry}{Entry}
\newcommand{\referxy}[4]{\bibitem{kn:#1}{#2,}~\emph{#3,}~{#4.}}	
\newcommand{\cito}[1]{\cite{kn:#1}}	
\newcommand{\citu}[2]{\cite[#2]{kn:#1}}

\begin{document}
\title{$q$-Analogues of $\pi$-Series by Applying
Carlitz Inversions to $q$-Pfaff-Saalsch\"utz Theorem}
\author{Xiaojing Chen and Wenchang Chu}
\address{School of Statistics\newline
	Qufu Normal University\newline
	Qufu (Shandong), P.~R.~China}
\email{upcxjchen@outlook.com}
\address{Department of Mathematics and Physics\newline
     	University of Salento,~P.~O.~Box~193 \newline
      	73100 Lecce, ~ Italy}
\email{chu.wenchang@unisalento.it}
\thanks{Correspondence: chu.wenchang@unisalento.it
	and upcxjchen@outlook.com}


\begin{abstract}
By applying multiplicate forms of the Carlitz inverse series
relations to the $q$-Pfaff-Saalsch\ütz summation theorem, we
establish twenty five nonterminating $q$-series identities 
with several of them serving as $q$-analogues of infinite 
series expressions for $\pi$ and $1/\pi$, including some 
typical ones discovered by Ramanujan (1914) and Guillera.
\end{abstract}

\maketitle\thispagestyle{empty}
\markboth{Xiaojing Chen and Wenchang Chu}{$q$-Analogues of $\pi$-Series}

\section{Introduction and Motivation}
Let $\mb{N}$ and $\mb{N}_0$ be the sets of natural numbers
and nonnegative integers, respectively. For an indeterminate
$x$, the Pochhammer symbol is defined by
\[(x)_0\equiv1\quad\text{and}\quad
(x)_n=x(x+1)\cdots(x+n-1)
\quad\text{for}\quad n\in\mb{N}\]
with the following shortened multiparameter notation
\[\hyp{cccc}{\alp,\bet,\cdots,\gam}{A,B,\cdots,C}_n
=\frac{(\alp)_n(\bet)_n\cdots(\gam)_n}{(A)_n(B)_n\cdots(C)_n}.\]
Analogously, the rising and falling $q$-factorials with the 
base $q$ are given by $(x;q)_0=\ang{x;q}_0\equiv1$ and 
\[\mult{l}{(x;q)_n=(1-x)(1-qx)\cdots(1-q^{n-1}x)\\[2mm]
\ang{x;q}_n=(1-x)(1-q^{-1}x)\cdots(1-q^{1-n}x)}
\Bigg\}\quad\text{for}\quad n\in\mb{N}.\]
Then the Gaussian binomial coefficient can be expressed as
\[\binq{m}{n}=\frac{(q;q)_m}{(q;q)_n(q;q)_{m-n}}
=\frac{(q^{m-n+1};q)_n}{(q;q)_n}
\quad\text{where}\quad m,\:n\in\mb{N}.\]
When $|q|<1$, the infinite product $(x;q)_{\infty}$ is well defined.
We have hence the $q$-gamma function~\citu{rahman}{\S1.10}
\[\Gam_q(x)\:=\:(1-q)^{1-x}
\frac{(q;q)_\infty}{(q^x;q)_\infty}
\quad\text{and}\quad
\lim_{q\to1^{-}}\Gam_q(x)\:=\:\Gam(x).\]
For the sake of brevity, the product and quotient of the $q$-shifted
factorials will be abbreviated respectively to
\bnm
\hh{\alp,\bet,\cdots,\gam}{q}{n}\quab
&=&\ff{\alp}{q}{n}\ff{\bet}{q}{n}\cdots\ff{\gam}{q}{n},\\
\hyq{cccc}{q}{\alp,\bet,\cdots,\gam}
          {A,B,\cdots,C}_n
&=&\frac{\ff{\alp}{q}{n}\ff{\bet}{q}{n}\cdots\ff{\gam}{q}{n}}
     {\ff{A}{q}{n}\ff{B}{q}{n}\cdots\ff{C}{q}{n}}.
\enm
According to Bailey~\cito{bailey} and Gasper--Rahman~\cito{rahman},
the $q$-series is defined by
\bnm
{_{1+\ell}\phi_\ell}
\hyq{rccc}{q;z}{a_0,a_1,\cdots,a_\ell}
                {b_1,\cdots,b_\ell}
=\:\:\sum_{n=0}^{\infty}
\hyq{cccc}{q}{a_0,a_1,\cdots,a_\ell}
     {q,b_1,\cdots,b_\ell}_n\:z^n.
\enm
This series is well defined when none of the denominator parameters
has the form $q^{-m}$ with $m\in\mb{N}_0$. If one of the numerator
parameters has the form $q^{-m}$ with $m\in\mb{N}_0$, the series
is terminating (in that case, it is a polynomial of $z$). Otherwise,
the series is said nonterminating, where we assume that $0<|q|<1$.

As the $q$-analogues of the Gould--Hsu~\cito{hsu} inversions,
Carlitz~\cito{carlitz73} found, in 1973, a well--known pair
of inverse series relations, which can be reproduced
as follows. Let $\{a_k,b_k\}_{k\ge0}$ be two sequences
such that the $\vph$-polynomials defined by
\[\label{psifunc}\vph(x;0)\equiv1
\quad\text{and}\quad
\vph(x;n)=\prod_{k=0}^{n-1}(a_k+xb_k)
\quad\text{for}\quad n=1,2,\cdots\]
differ from zero at $x=q^{-m}$ for $m\in\mb{N}_0$.
Then the following inverse relations hold
\pmq{
f(n)=&\label{carlitz-f}
\sum_{k=0}^n(-1)^k
\binq{n}{k}\vph(q^{-k};n)g(k),\\
g(n)=&\label{carlitz-g}
\sum_{k=0}^n(-1)^k
\binq{n}{k}q^{\binm{n-k}{2}}
\frac{a_k+q^{-k}b_k}{\vph(q^{-n};k+1)}f(k).}	
Alternatively, if the $\vph$-polynomials differ from zero
at $x=q^{m}$ for $m\in\mb{N}_0$, Carlitz deduced, under
the base change $q\to q^{-1}$, another equivalent pair
\pmq{
f(n)=&\label{carlitz+f}
\sum_{k=0}^n(-1)^k
\binq{n}{k}q^{\binm{n-k}{2}}\vph(q^{k};n)g(k),\\
g(n)=&\label{carlitz+g}
\sum_{k=0}^n(-1)^k\binq{n}{k}
\frac{a_k+q^{k}b_k}{\vph(q^{n};k+1)}f(k).}
These inverse series relations have been shown by 
and Chu~\cite{kn:chu94d,kn:chu95a} to be very useful 
in proving terminating $q$-series identities.
Among numerous $q$-series identities, the
$q$-Pfaff--Saalsch\"utz theorem (cf.~\citu{rahman}{II-12}
for the terminating balanced series is fundamental:
\pq{\label{pfaff-q}{_3\phi_2}\hyq{ccc}{q;q}
{q^{-n},a,b}{c,q^{1-n}ab/c}
=\hyq{cc}{q}{c/a,c/b}{c,c/ab}_n.}

As a warm--up, we illustrate how to derive the $q$-Dougall
sum by making use of Carlitz' inversions. Observe that
\eqref{pfaff-q} is equivalent to
\[{_3\phi_2}\hyq{ccc}{q;q}
{q^{-n},q^na,qa/bd}{qa/b,qa/d}
=\Big(\frac{qa}{bd}\Big)^n
\hyq{cc}{q}{b,~d}{qa/b,qa/d}_n\]
which can be rewritten as a $q$-binomial sum
\[\sum_{k=0}^n(-1)^k\binq{n}{k}
q^{\binm{n-k}2}(q^ka;q)_n
\hyq{c}{q}{a,qa/bd}{qa/b,qa/d}_k
=\Big(\frac{qa}{bd}\Big)^n
\hyq{c}{q}{a,~b,~d}{qa/b,qa/d}_n
q^{\binm{n}2}.\]
This matches exactly \eqref{carlitz+f}
under the specifications
\pnq{f(n)&=\Big(\frac{qa}{bd}\Big)^n
\hyq{c}{q}{a,~b,~d}{qa/b,qa/d}_nq^{\binm{n}2},\\
g(k)&=\hyq{c}{q}{a,~qa/bd}{qa/b,qa/d}_k
\quad\text{and}\quad
\vph(x;n)=(ax;q)_n.}
Then the dual relation corresponding to \eqref{carlitz+g}
reads as
\[\sum_{k=0}^n(-1)^k\binq{n}{k}
\frac{1-q^{2k}a}{(q^na;q)_{k+1}}
\Big(\frac{qa}{bd}\Big)^k
\hyq{c}{q}{a,~b,~d}{qa/b,qa/d}_kq^{\binm{k}2}
=\hyq{c}{q}{a,~qa/bd}{qa/b,qa/d}_n.\]
This is equivalent to the $q$-Dougall sum (cf.~\citu{rahman}{II-21}):
\[{_6\phi_5}
\hyq{rrccc}{q;\frac{q^{n+1}a}{bd}}
{a,~q\sqrt{a},&-q\sqrt{a},&b,~\qquad~d,~\qquad~q^{-n}}
{\sqrt{a},&-\sqrt{a},&qa/b,qa/d,q^{n+1}a}
=\hyq{c}{q}{qa,qa/bd}{qa/b,qa/d}_n.\]

For $a=b=d=q^{1/2}$, the limiting case $n\to\infty$
of the last formula becomes
\[\frac1{\Gam^2_q(\frac12)}=\sum_{k=0}^{\infty}
(-1)^k\frac{(q^{1/2};q)_k^3}{(q;q)_k^3}
\frac{1-q^{2k+\frac12}}{1-q}
q^{\frac{k^2}2}\]
which  reduces, for $q\to1^{-}$, to the following infinite
series expression for $\pi$
\[\frac{2}{\pi}=\sum_{k=0}^{\infty}
(-1)^k\frac{(\frac12)^3_k}{(1)^3_k}\big\{1+4k\big\}\]
as recorded in one of Ramanujan's letters to Hardy~\cito{hardy}.
More difficult formulae for $1/\pi$ were subsequently
discovered by Ramanujan~\citu{ramanujan}{1914}, where
17 similar series representations were announced.
Three of them are reproduced as follows:
\pnq{\frac{4}{\pi}
=&\sum_{k=0}^{\infty}
\hyp{ccc}
{\frac12,\frac12,\frac12}{\rule[1mm]{0mm}{3mm}1,\:1,\:1}_k
\frac{1+6k}{4^k}.\\
\frac{8}{\pi}
=&\sum_{k=0}^{\infty}
\hyp{ccc}
{\frac12,\frac14,\frac34}
{\rule[1mm]{0mm}{3mm}1,\:1,\:1}_k
\frac{3+20k}{(-4)^k}.\\
\frac{16}{\pi}\!
=&\sum_{k=0}^{\infty}
\hyp{ccc}
{\frac12,\frac12,\frac12}{\rule[1mm]{0mm}{3mm}1,\:1,\:1}_k
\frac{5+42k}{64^k}.}
For their proofs and recent developments, the reader can consult
the papers by Baruah-Berndt-Chan~\cito{berndt-bc},
Guillera~\cite{kn:gj03em,kn:gj06rj,kn:gj08rj} and  
Chu~\emph{et al}~\cite{kn:chu11a,kn:chu14f}.

Recently, there has been a growing interest in finding
$q$-analogues of Ramanujan--like series
(cf.~\cite{kn:chu18e,kn:chu20RJ,
	kn:guo18jdea,kn:guo18itsf,kn:guo19rj,kn:gj18jdea}).
Following the procedure just described, the aim of this paper
is to show systematically $q$-analogues of $\pi$-related series
by applying the multiplicate form of Carlitz inverse series
relations to the $q$-Pfaff--Saalsch\"utz summation theorem.
In the next section, we shall derive, by employing the duplicate
inversions, twenty $q$-series identities including $q$-analogues
of the afore displayed three series of Ramanujan. Then in section~3,
the triplicate inversions will be utilized to establish five $q$-series
identities. By applying the bisection series method to two resulting
series, $q$-analogues are established also for the following two
remarkable series discovered by Guillera~\cite{kn:gj03em,kn:gj06rj}:
\pnq{\frac{2\sqrt2}{\pi}~&
=\sum_{k=0}^\infty~\Big(\frac{-1}8\Big)^k
\hyp{ccc}{\frac12,\frac12,\frac12}
{1,\:1,\:1\rule[2mm]{0mm}{2mm}}_k
\big\{1+6k\big\}.\\
\frac{32\sqrt2}{\pi}
&=\sum_{k=0}^{\infty}\Big(\frac{-3}{8}\Big)^{3k}
\hyp{ccccc}{\frac12,\frac16,\frac56\\[-3mm]}
{1,\:1,\:1}_k\big\{15+154k\big\}.}

\section{Duplicate Inverse Series Relations}

Denote by $\pav{x}$ the integer part for a real number $x$.
Then for all the $n\in\mb{N}_0$, there holds the equality
\[\boxed{n=\pavv{\tfrac{n}2}+\pavv{\tfrac{1+n}2}}.\]
According to this partition, we shall reformulate \eqref{pfaff-q}
in three different manners. Their dual relations will lead us
to $q$-series counterparts for several remarkable infinite
series expressions of $\pi$ and $1/\pi$.

\subsection{} \
According to the $q$-Pfaff--Saalsch\"utz formula \eqref{pfaff-q},
it is not hard to verify that
\[{_3\phi_2}\hyq{ccc}{q;q}{q^{-n},~a,~c}
{q^{-\pav{\frac{n}2}}ae,q^{1-\pav{\frac{n+1}2}}c/e}
=\hyq{c}{q}{q^{-\pav{\frac{n}2}}e,q^{-\pav{\frac{n}2}}ae/c}
{q^{-\pav{\frac{n}2}}ae,q^{-\pav{\frac{n}2}}e/c}_n\]
which is equivalent to the binomial sum
\pnq{
&\sum_{k=0}^n
(-1)^k\binq{n}{k}
(q^{1-k}/ae;q)_{\pav{\frac{n}{2}}}
(q^{-k}e/c;q)_{\pav{\frac{n+1}{2}}}
\hyq{c}{q}{a,\:c}{ae,qc/e}_kq^{\binm{k+1}2}\\
&\:=\hyq{c}{q}{e,ae/c}{ae}_{\pav{\frac{n+1}{2}}}
\hyq{c}{q}{q/e,qc/ae}{qc/e}_{\pav{\frac{n}{2}}}.}
Observing that this equation matches exactly
to \eqref{carlitz-f} specified by
\pnq{
&f(k)=\hyq{c}{q}{e,ae/c}{ae}_{\pav{\frac{k+1}{2}}}
\hyq{c}{q}{q/e,qc/ae}{qc/e}_{\pav{\frac{k}{2}}},\\
&g(k)=\hyq{c}{q}{a,\:c}{ae,qc/e}_k
q^{\binm{k+1}2},\\[2mm]
&\vph(x;n)=(qx/ae;q)_{\pav{\frac{n}{2}}}
(ex/c;q)_{\pav{\frac{n+1}{2}}};}
we may state the dual relation corresponding
to \eqref{carlitz-g} as the proposition.
\begin{prop}[Terminating reciprocal relation]\label{pp=a}
\pnq{\hyq{c}{q}{a,\:c}{ae,qc/e}_n
=\sum_{k\ge0}\binq{n}{2k}
\frac{(1-q^{-k}e/c)q^{(1+2k)(k-n)}}
     {(q^{1-n}/ae;q)_{k}(q^{-n}e/c;q)_{k+1}}
\hyq{c}{q}{e,ae/c}{ae}_{k}
\hyq{c}{q}{q/e,qc/ae}{qc/e}_{k}&\\
-\sum_{k\ge0}\binq{n}{2k+1}
\frac{(1-q^{-k}/ae)q^{(1+k)(1+2k-2n)}}
     {(q^{1-n}/ae;q)_{k+1}(q^{-n}e/c;q)_{k+1}}
\hyq{c}{q}{e,ae/c}{ae}_{k+1}
\hyq{c}{q}{q/e,qc/ae}{qc/e}_{k}&.}
\end{prop}
The two sums just displayed are, in fact, balanced
$_8\phi_7$-series, which do not admit closed forms.
However their combination does have a closed form.
That is the reason why we call the last relation reciprocal.

Letting $n\to\infty$ in Proposition~\ref{pp=a} and then
applying the Weierstrass $M$-test on uniformly convergent
series (cf. Stromberg~\citu{karl}{\S3.106}), we get the
limiting relation:
\pnq{\hyq{c}{q}{a,\:c}{ae,qc/e}_{\infty}
=\sum_{k\ge0}
\frac{1-q^kc/e}{(q;q)_{2k}}
\hyq{c}{q}{e,ae/c}{ae}_{k}
\hyq{c}{q}{q/e,qc/ae}{qc/e}_{k}
q^{k^2-k}(ac)^k&\\
+\frac{c}{e}\sum_{k\ge0}
\frac{1-ae/q}{(q;q)_{2k+1}}
\hyq{c}{q}{e,ae/c}{ae/q}_{k+1}
\hyq{c}{q}{q/e,qc/ae}{qc/e}_{k}
q^{k^2}(ac)^k&.}

By unifying the two sums together, we find the following theorem.
\begin{thm}[Nonterminating series identity]\label{thm=a}
\pnq{\hyq{c}{q}{a,\:c}{ae,c/e}_{\infty}
=\sum_{k=0}^{\infty}
&\frac{(ac)^k}{(q;q)_{2k}}
\hyq{c}{q}{e,ae/c}{ae}_{k}
\hyq{c}{q}{q/e,qc/ae}{c/e}_{k}
q^{k^2-k}\\
\times~&\bigg\{1+q^k\frac{c(1-q^ke)(1-q^kae/c)}
	{e(1-q^{1+2k})(1-q^kc/e)}\bigg\}.}
\end{thm}

We highlight two important corollaries about product of reciprocal
$q$-gamma functions. Their limiting case $q\to1^{-}$ will yield infinite
series for $\pi$ and $1/\pi$.
\begin{corl}[$a=q^{\lam}$ and $c=e=q^{1-\lam}$
            in Theorem~\ref{thm=a}]\label{cc-a}
\pnq{\frac{1}{\Gam_q(1+\lam)\Gam_q(2-\lam)}
=\sum_{k=0}^{\infty}&
\frac{\hh{q^{\lam},q^{1+\lam},q^{1-\lam},q^{2-\lam}}{q}{k}}
     {(q;q)^2_k(q^2;q)_{2k}}q^{k^2+k}\\
\times~&\bigg\{1-\frac{(1-q^{-k}) (1-q^{1+2k})}
{(1-q^{\lambda+k})(1-q^{1-\lambda+k})}\bigg\}.}
\end{corl}

\begin{corl}[$a=c=q$ and $e=q^{\lam}$ in Theorem~\ref{thm=a}]\label{cc+a}
\[\Gam_q(\lam)\Gam_q(1-\lam)
=\sum_{k=0}^{\infty}q^{k^2+k}
\frac{(q^{\lam};q)_k(q^{1-\lam};q)_k}{(q^2;q)_{2k}}
\bigg\{\frac{1-q^{1+2k}}{1-q^{\lam+k}}
-\frac{1-q^{\lam+k}}{1-q^{\lam-1-k}}\bigg\}.\]
\end{corl}

By properly choosing special values of $a,c$ and $e$, we find ten
interesting $q$-series identities, that correspond to the classical
series with convergence rate ``$\frac14$".

\textbf{A1}. \
For the series discovered by Ramanujan~\cito{ramanujan}
\[\frac{4}{\pi}=\sum_{k=0}^\infty
\hyp{ccc}{\frac12,\frac12,\frac12}{1,\:1,\:1}_k
\frac{1+6k}{4^k},\]
we recover, by letting $\lam=1/2$ in Corollary~\ref{cc-a}, 
the following $q$-analogue (cf.~Chen--Chu~\citu{chu20RJ}{Example~38} 
and Guo~\citu{guo19rj}{Equation~1.6})
\[\frac{1}{\Gam^2_q(\frac12)}
=\sum_{k=0}^{\infty}q^{k^2}
\frac{(q^{1/2};q)^4_k}{(q;q)^2_k(q;q)_{2k}}
\frac{1+q^{k+1/2}-2q^{2k+1/2}}
     {(1-q)(1+q^{k+1/2})}.\]
A different, but simpler $q$-analogue can be 
found in Guo--Liu~\citu{guo18jdea}{Equation~3} 
and Chen--Chu~\citu{chu20RJ}{Example~4}:
\[\sum_{k=0}^{\infty}
\frac{1-q^{6k+1}}{1-q^4}
\frac{(q;q^2)^2_k(q^2;q^4)_k}
     {(q^4;q^4)^3_k}
q^{k^2}
=\frac{1}{\Gam^2_{q^4}(\frac12)}.\]

\textbf{A2}. \
For $\lam=1/3$, we get, from Corollary~\ref{cc-a}, the following series
\[\frac{1}{\Gam_q(\frac43)\Gam_q(\frac53)}
=\sum_{k=0}^{\infty}q^{k^2+k}
\frac{\hh{q^{1/3},q^{2/3},q^{4/3},q^{5/3}}{q}{k}}{(q;q)^2_k(q^2;q)_{2k}}
\bigg\{1-\frac{(1-q^{-k})(1-q^{2 k+1})}
{(1-q^{k+\frac{1}{3}})(1-q^{k+\frac{2}{3}})}\bigg\}\]
which gives a $q$-analogue of the series
\[\frac{9\sqrt3}{2\pi}=\sum_{k=0}^\infty
\hyp{cccc}{\frac13,\frac13,\frac23,\frac23}
    {1,1,1,\:\frac32\rule[2mm]{0mm}{2mm}}_k
\frac{2+18k+27k^2}{4^k}.\]

\textbf{A3}. \
For $\lam=1/4$, we have, from Corollary~\ref{cc-a}, the following series
due to Guo and Zudilin~\citu{guo18itsf}{Equation~1.6}
\[\frac{1}{\Gam_q(\frac54)\Gam_q(\frac74)}
=\sum_{k=0}^{\infty}q^{k^2+k}
\frac{\hh{q^{\frac14},q^{\frac34},q^{\frac54},q^{\frac74}}{q}{k}}{(q;q)^2_k(q^2;q)_{2k}}
\bigg\{1-\frac{(1-q^{-k})(1-q^{2k+1})}{(1-q^{k+\frac14})(1-q^{k+\frac34})}\bigg\}\]
which offers a $q$-analogue of the series
\[\frac{8\sqrt2}{\pi}=\sum_{k=0}^\infty
\hyp{cccc}{\frac14,\frac14,\frac34,\frac34}
    {1,1,1,\:\frac32\rule[2mm]{0mm}{2mm}}_k
\frac{3+32k+48k^2}{4^k}.\]

\textbf{A4}. \
For $\lam=1/6$, we find, from Corollary~\ref{cc-a}, the following series
\[\frac{1}{\Gam_q(\frac76)\Gam_q(\frac{11}6)}
=\sum_{k=0}^{\infty}q^{k^2+k}
\frac{\hh{q^{1/6},q^{5/6},q^{7/6},q^{11/6}}{q}{k}}{(q;q)^2_k(q^2;q)_{2k}}
\bigg\{1-\frac{(1-q^{-k})(1-q^{2 k+1})}
{(1-q^{k+\frac{1}{6}})(1-q^{k+\frac{5}{6}})}\bigg\}\]
which provides a $q$-analogue of the series
\[\frac{18}{\pi}=\sum_{k=0}^\infty
\hyp{cccc}{\frac16,\frac16,\frac56,\frac56}
    {1,1,1,\frac32\rule[2mm]{0mm}{2mm}}_k
\frac{5+72k+108k^2}{4^k}.\]

\textbf{A5}. \
Letting $\lam=1/2$ in Corollary~\ref{cc+a}, we get the following series
\[\Gam^2_q\big(\frac12\big)
=\sum_{k=0}^{\infty}q^{k^2+k}
\frac{(q^{1/2};q)^2_k}{(q^2;q)_{2k}}
(1+2q^{k+1/2})\]
which is a $q$-analogue of the series
\[\frac{\pi}{3}=\sum_{k=0}^\infty
\hyp{ccc}{\frac12,\frac12\\[-3mm]}{1,\frac32}_k
\bigg(\frac14\bigg)^k.\]

\textbf{A6}. \
Letting $\lam=1/3$ in Corollary~\ref{cc+a}, we deduce the following series
\[\Gam_q\bigg(\frac13\bigg)\Gam_q\bigg(\frac23\bigg)
=\sum_{k=0}^{\infty}q^{k^2+k}
\frac{(q^{1/3};q)_k(q^{2/3};q)_k}{(q^2;q)_{2k}}
\bigg\{\frac{1-q^{1+2k}}{1-q^{k+\frac13}}
-\frac{1-q^{k+\frac13}}{1-q^{-k-\frac23}}\bigg\}\]
which gives a $q$-analogue of the series
\[\frac{4\pi}{\sqrt{3}}=\sum_{k=0}^\infty
\hyp{cccc}{\frac13,\frac13,\frac23,\frac23}
    {1,\frac32,\frac43,\frac53\rule[2mm]{0mm}{2mm}}_k
\frac{7+27k+27k^2}{4^k}.\]

\textbf{A7}. \
Letting $\lam=1/6$ in Corollary~\ref{cc+a}, we obtain the following series
\[\Gam_q\bigg(\frac16\bigg)\Gam_q\bigg(\frac56\bigg)
=\sum_{k=0}^{\infty}q^{k^2+k}
\frac{(q^{1/6};q)_k(q^{5/6};q)_k}{(q^2;q)_{2k}}
\bigg\{\frac{1-q^{1+2k}}{1-q^{k+\frac16}}
-\frac{1-q^{k+\frac16}}{1-q^{-k-\frac56}}\bigg\}\]
which results in a $q$-analogue of the series
\[{10\pi}=\sum_{k=0}^\infty
\hyp{cccc}{\frac16,\frac16,\frac56,\frac56}
    {1,\frac32,\frac76,\frac{11}6\rule[2mm]{0mm}{2mm}}_k
\frac{31+108k+108k^2}{4^k}.\]

\textbf{A8}. \
By specifying $a=q,~c=q^{2/3}$ and $e=q^{1/3}$
in Theorem~\ref{thm=a}, we find
\[\frac{\Gam_q^2(\frac13)}{\Gam_q(\frac23)}
=\sum_{k=0}^{\infty}q^{k^2+\frac{2k}3}
\frac{(q^{1/3};q)_k(q^{2/3};q)^2_k}{(q^{4/3};q)_k(q^2;q)_{2k}}
\frac{1+q^{k+\frac{1}{3}}-2q^{2k+1}}{1-q^\frac13}\]
which corresponds to the classical series
\[\frac{\sqrt{3}\:\Gam^3(\frac13)}{2\pi}
=\sum_{k=0}^\infty
\hyp{cccc}{\frac13,\frac23,\frac23}
    {1,\frac32,\frac43\rule[2mm]{0mm}{2mm}}_k
\frac{5+9k}{4^k}.\]

\textbf{A9}. \
By specifying $a=c=q^{1/4}$ and $e=q^{1/2}$
in Theorem~\ref{thm=a}, we have
\[\frac{\Gam_q^2(\frac34)}{\Gam_q^2(\frac14)}
=\sum_{k=0}^{\infty}q^{k(k-\frac12)}
\frac{(q^{1/2};q)^3_k(q^{3/2};q)_k}{(q^{3/4};q)^2_k(q^2;q)_{2k}}
\frac{1+q^{k+\frac12}-2q^{2k+\frac14}}{(1-q)(1+q^\frac12)}\]
which corresponds to the classical series
\[\frac{2\Gam^2(\frac34)}{3\Gam^2(\frac14)}
=\sum_{k=0}^\infty
\hyp{cccc}{\frac12,\frac12,\frac12}
    {1,\frac34,\frac34\rule[2mm]{0mm}{2mm}}_k
\frac{k}{4^k}
~\Longleftrightarrow~
\frac{12\Gam^2(\frac34)}{\Gam^2(\frac14)}
=\sum_{k=0}^\infty
\hyp{cccc}{\frac32,\frac32,\frac32}
    {1,\frac74,\frac74\rule[2mm]{0mm}{2mm}}_k
\frac{1}{4^k}.\]

\textbf{A10}. \
By specifying $a=c=q^{3/4}$ and $e=q^{1/2}$
in Theorem~\ref{thm=a}, we find
\[\frac{\Gam_q^2(\frac14)}{\Gam_q^2(\frac34)}
=\sum_{k=0}^{\infty}q^{k(k+\frac12)}
\frac{(q^{1/2};q)^3_k(q^{3/2};q)_k}{(q^{5/4};q)^2_k(q^2;q)_{2k}}
\frac{(1+q^\frac14)(1+q^{k+\frac12}-2q^{2k+\frac34})}{1-q^\frac14}\]
which corresponds to the classical series
\[\frac{\Gam^2(\frac14)}{8\Gam^2(\frac34)}
=\sum_{k=0}^\infty
\hyp{cccc}{\frac12,\frac12,\frac12}
    {1,\frac54,\frac54\rule[2mm]{0mm}{2mm}}_k
\frac{1+3k}{4^k}.\]

\subsection{} \
According to \eqref{pfaff-q}, it is routine to check that
\[{_3\phi_2}\hyq{ccc}{q;q}{q^{-n},~q^{\pav{\frac{n}2}}a,~c}
{ae,q^{1-\pav{\frac{n+1}2}}c/e}
=\hyq{c}{q}{q^{-\pav{\frac{n}2}}e,ae/c}
{q^{-\pav{\frac{n}2}}e/c,ae}_n.\]
By making use of the factorial expression
\[(q^{-k}y;q)_{\pav{\frac{n+1}2}}q^{\pav{\frac{n+1}2}k}
=\ang{q^k/y;q}_{\pav{\frac{n+1}2}}
(-y)^{\pav{\frac{n+1}2}}
q^{\binm{\pav{\frac{n+1}2}}{2}},\]
we can reformulate the last equality as the $q$-binomial identity:
\pnq{
&\sum_{k=0}^n
(-1)^k\binq{n}{k}q^{\binm{n-k}2}
(q^{k}a;q)_{\pav{\frac{n}{2}}}
\ang{q^{k}c/e;q}_{\pav{\frac{n+1}{2}}}
\hyq{c}{q}{a,\:c}{ae,qc/e}_k\\
&\:=(-1)^{\pav{\frac{n+1}{2}}}
q^{\binm{n}2-\binm{\pav{\frac{n+1}2}}{2}}c^n
\frac{(e;q)_{\pav{\frac{n+1}2}}}{e^{\pav{\frac{n+1}{2}}}}
\hyq{c}{q}{ae/c}{ae}_{n}
\hyq{c}{q}{q/e,a}{qc/e}_{\pav{\frac{n}{2}}}.}
Since the last equation matches exactly
to \eqref{carlitz+f} specified by
\pnq{
f(k)&=(-1)^{\pav{\frac{k+1}{2}}}
q^{\binm{k}2-\binm{\pav{\frac{k+1}2}}{2}}c^k
\frac{(e;q)_{\pav{\frac{k+1}2}}}{e^{\pav{\frac{k+1}{2}}}}
\hyq{c}{q}{ae/c}{ae}_{k}
\hyq{c}{q}{q/e,a}{qc/e}_{\pav{\frac{k}{2}}},\\
g(k)&=\hyq{c}{q}{a,\:c}{ae,qc/e}_k
\quad\text{and}\quad
\vph(x;n)=(ax;q)_{\pav{\frac{n}{2}}}
\ang{cx/e;q}_{\pav{\frac{n+1}{2}}};}
the dual relation corresponding to \eqref{carlitz+g}
is given in the proposition.
\begin{prop}[Terminating reciprocal relation]\label{pp=b}
\pnq{\hyq{c}{q}{a,\:c}{ae,qc/e}_n
=\sum_{k\ge0}q^{\frac{3k^2-k}2}\binq{n}{2k}
\frac{(1-q^{k}c/e)(-1)^kc^{2k}}
     {(q^{n}a;q)_{k}\ang{q^{n}c/e;q}_{k+1}}
\frac{(e;q)_{k}}{e^{k}}
\hyq{c}{q}{ae/c}{ae}_{2k}
\hyq{c}{q}{q/e,a}{qc/e}_{k}&\\
+\sum_{k\ge0}q^{\frac{3k^2+k}2}\binq{n}{2k+1}
\frac{(1-aq^{3k+1})(-1)^kc^{2k+1}}
     {(q^{n}a;q)_{k+1}\ang{q^{n}c/e;q}_{k+1}}
\frac{(e;q)_{k+1}}{e^{k+1}}
\hyq{c}{q}{ae/c}{ae}_{2k+1}
\hyq{c}{q}{q/e,a}{qc/e}_{k}&.}
\end{prop}
Both sums just displayed can be expressed as terminating
$q$-series, which do not have closed forms. However
their combination does have a closed form.

Letting $n\to\infty$ in Proposition~\ref{pp=b} and then
applying the Weierstrass $M$-test on uniformly convergent
series,
we get the
limiting relation:
\pnq{\hyq{c}{q}{a,\:c}{ae,qc/e}_{\infty}
=\sum_{k\ge0}
(-1)^kq^{\frac{3k^2-k}2}
\frac{1-q^kc/e}{(q;q)_{2k}}
\frac{c^{2k}(e;q)_{k}}{e^{k}}
\hyq{c}{q}{ae/c}{ae}_{2k}
\hyq{c}{q}{q/e,a}{qc/e}_{k}&\\
+\frac{c}{e}\sum_{k\ge0}
(-1)^kq^{\frac{3k^2+k}2}
\frac{1-aq^{3k+1}}{(q;q)_{2k+1}}
\frac{c^{2k}(e;q)_{k+1}}{e^{k}}
\hyq{c}{q}{ae/c}{ae}_{2k+1}
\hyq{c}{q}{q/e,a}{qc/e}_{k}&.}

By unifying the two sums together, we find the following theorem.
\begin{thm}[Nonterminating series identity]\label{thm=b}
\pnq{\hyq{c}{q}{a,\:c}{ae,c/e}_{\infty}
&=\sum_{k=0}^{\infty}
\frac{(-c^2/e)^k}{(q;q)_{2k}}
\frac{(ae/c;q)_{2k}}{(ae;q)_{2k}}
\hyq{c}{q}{a,e,q/e}{c/e}_{k}
q^{\frac{3k^2-k}2}\\
&\times\bigg\{1+q^k\frac{c(1-aq^{3k+1})(1-q^ke)(1-q^{2k}ae/c)}
	{e(1-q^{1+2k})(1-q^kc/e)(1-aeq^{2k})}\bigg\}.}
\end{thm}

Two implications are given below about product
of reciprocal $q$-gamma functions.
\begin{corl}[$a=q^{\lam}$ and $c=e=q^{1-\lam}$ in Theorem~\ref{thm=b}]\label{cc-b}
\pnq{\frac1{\Gam_q(1+\lam)\Gam_q(2-\lam)}
&=\sum_{k=0}^{\infty}(-1)^k
\frac{(q^{1+\lam};q)_{2k}}{(q^2;q)^2_{2k}}
\hyq{c}{q}{q^{\lam},q^{\lam},q^{2-\lam}}{q}_{k}
q^{\frac{k}2(3+3k-2\lam)}\\
&\times\frac{1-q^{1+\lam+3k}}{1-q}\bigg\{1+
\frac{q^{-k}(1-q^k)(1-q^{1+2k})(1-q^{1+2k})}
     {(1-q^{1-\lam+k})(1-q^{\lam+2k})(1-q^{1+\lam+3k})}\bigg\}.}
\end{corl}

\begin{corl}[$a=c=q$ and $e=q^{\lam}$ in Theorem~\ref{thm=b}]\label{cc+b}
\pnq{{\Gam_q(1+\lam)\Gam_q(1-\lam)}
&=\sum_{k=0}^{\infty}(-1)^k
\frac{\hh{q,q^{\lam}}{q}{k}}{(q;q)_{2k}}
\frac{(q^{\lam};q)_{2k}}{(q^{1+\lam};q)_{2k}}
q^{\frac{k}2(3+3k-2\lam)}\\
&\times\bigg\{1+
\frac{q^{1+k-\lam}(1-q^{2+3k})(1-q^{\lam+k})(1-q^{\lam+2k})}
     {(1-q^{1+2k})(1-q^{1-\lam+k})(1-q^{1+\lam+2k})}\bigg\}.}
\end{corl}

Five $q$-series as well as their counterparts
of classical series are exemplified as follows.

\textbf{B1}. \
For Ramanujan's series~\cito{ramanujan}
 \[\frac{8}{\pi}=\sum_{k=0}^\infty
\Big(\frac{-1}{4}\Big)^k
\hyp{ccc}{\frac12,\frac14,\frac34}{1,\,1,\,1}_k
\big\{3+20k\big\},\]
we recover, by letting $\lam=1/2$ in Corollary~\ref{cc-b},
the following $q$-analogue (cf.~Chen and Chu~\citu{chu20RJ}{Example~39})
\pnq{\frac{1}{\Gam^2_q(\frac12)}
=&\sum_{k=0}^{\infty}(-1)^k
\frac{(q^{1/2};q)^3_k(q^{1/2};q)_{2k}}{(q;q)_k(q;q)^2_{2k}}
q^{3k^2/2}\\
\times&\nnm\bigg\{
\frac{(1+q^{k+1/2})^2(1-q^{3k+1/2})-q^{2k+1/2}(1-q^{2k+1/2})}
     {(1-q)(1+q^{k+1/2})^2}\bigg\}.}

Guo and Zudilin~\citu{guo18itsf}{Equation~1.4} derived,
by means of the \emph{WZ} machinery, another $q$-analogue
\pnq{\frac{1}{\Gam^2_q(\frac12)}
=&\sum_{k=0}^{\infty}(-1)^k
\frac{(q^{1/2};q)^2_k(q^{1/4};q^{1/2})_{2k}}{(q;q)^2_k(q;q)_{2k}}
q^{k^2/2}\\
\times&\nnm\bigg\{
\frac{1-q^{2k+1/4}}{1-q}+
\frac{q^{k+1/4}(1-q^{k+1/4})}{(1-q)(1+q^{k+1/2})}\bigg\}.}

This is another example (apart from \textbf{A1}) that there 
may exist different $q$-analogues for some classical series.

\textbf{B2}. \
For $\lam=1/2$, we get, from Corollary~\ref{cc+b}, the series
\pnq{\Gam_q\bigg(\frac12\bigg)\Gam_q\bigg(\frac32\bigg)
=\sum_{k=0}^{\infty}&(-1)^kq^{\frac{3 k^2}{2}+k}
\frac{(q;q)_k(q^{1/2};q)_k(q^{1/2};q)_{2k}}{(q^{3/2};q)_{2k}(q;q)_{2k}}\\
\times~&\bigg\{1+
\frac{q^{k+\frac12}(1-q^{3k+2})(1-q^{2k+\frac12})}
     {(1-q^{2k+1})(1-q^{2k+\frac32})}\bigg\}}
which can also be obtained from Chu~\citu{chu18e}{Proposition~14: $x=y^2=q$}.
The above series gives, in turn, a $q$-analogue of the classical series
\[\frac{3\pi}{2}
=\sum_{k=0}^\infty\Big(\frac{-1}{4}\Big)^k
\hyp{ccc}{\frac12,\frac14,\frac34}
    {\frac32,\frac54,\frac74\rule[2mm]{0mm}{2mm}}_k
\big\{5+21k+20k^2\big\}.\]

\textbf{B3}. \
For $\lam=1/3$, we have, from Corollary~\ref{cc+b}, the series
\pnq{\Gam_q
\bigg(\frac13\bigg)\Gam_q\bigg(\frac23\bigg)
&=\sum_{k=0}^{\infty}(-1)^{k+1}
\frac{(q;q)_{k+1}(q^{1/3};q)_k(q^{4/3};q)_{2k}}
	{(q^2;q)_{2k}(q^{4/3};q)_{2k+1}}
q^{\frac{3k^2}2+\frac{19k}6+1}\\
&\times
\bigg\{1+\frac{(1+q^{k+\frac{2}{3}})(1-q^{-2 k-1}) (1-q^{3 k+1})}
   {(1-q^{k+1}) (1-q^{2 k+\frac{1}{3}})}\bigg\}}
which offers a $q$-analogue of the classical series
\[\frac{8\pi}{3\sqrt3}
=\sum_{k=0}^\infty\Big(\frac{-1}{4}\Big)^k
\hyp{ccc}{\frac13,\frac23,\frac16}
    {\frac32,\frac53,\frac76\rule[2mm]{0mm}{2mm}}_k
\big\{5+23k+30k^2\big\}.\]

\textbf{B4}. \
For $\lam=2/3$, we obtain, from Corollary~\ref{cc+b}, the series
\pnq{\Gam_q
\bigg(\frac13\bigg)\Gam_q\bigg(\frac53\bigg)
&=\sum_{k=0}^{\infty}(-1)^{k}
\frac{(q;q)_{k}(q^{2/3};q)_k(q^{2/3};q)_{2k}}
	{(q;q)_{2k}(q^{5/3};q)_{2k}}
q^{\frac{3k^2}2+\frac{5k}6}\\
&\times
\bigg\{1+\frac{q^{k+\frac13}(1+q^{k+\frac{1}{3}})(1-q^{k+\frac23}) (1-q^{3 k+2})}
   {(1-q^{2k+1}) (1-q^{2 k+\frac{5}{3}})}\bigg\}.}
which provides a $q$-analogue of the classical series
\[\frac{20\pi}{3\sqrt3}
=\sum_{k=0}^\infty\Big(\frac{-1}{4}\Big)^k
\hyp{ccc}{\frac13,\frac23,\:\frac56}
    {\frac32,\frac43,\frac{11}6\rule[2mm]{0mm}{2mm}}_k
\big\{13+40k+30k^2\big\}.\]

\textbf{B5}. \
In addition, by specifying $a=q^{2/3},~c=q^{1/3}$ and $e=q^{1/6}$ in Theorem~\ref{thm=b},
we find the following strange looking identity	
\pnq{\frac{\Gam_q(\frac16)\Gam_q(\frac56)}
{\Gam_q(\frac13)\Gam_q(\frac23)}
&=\sum_{k=0}^{\infty}(-1)^k
\frac{\hh{q^{\frac23},q^{\frac56}}{q}{k}}{(q;q)_{2k}}
\frac{(q^{\frac12};q)_{2k}}{(q^{\frac56};q)_{2k}}
q^{\frac{3k^2}2}\\
&\times\bigg\{1+q^{k+\frac16}
\frac{(1-q^{\frac12+2k})(1-q^{\frac53+3k})}
     {(1-q^{1+2k})(1-q^{\frac56+2k})}\bigg\}}
which turns out to be a $q$-analogue of the classical series
\[{5\sqrt3}
=\sum_{k=0}^\infty\Big(\frac{-1}{4}\Big)^k
\hyp{cccc}{\frac23,\:\frac14,\frac34,\:\frac56}
    {1,\frac32,\frac{11}{12},\frac{17}{12}\rule[2mm]{0mm}{2mm}}_k
\big\{10+51k+60k^2\big\}.\]

\subsection{} \
According to \eqref{pfaff-q}, it is not difficult to show that
\[{_3\phi_2}\hyq{ccc}{q;q}
{q^{-n},q^{\pav{\frac{n}2}}a,q^{\pav{\frac{n+1}2}}c}
{ae,\quad qc/e}
=\hyq{c}{q}{q^{-\pav{\frac{n}2}}e,q^{-\pav{\frac{n+1}2}}ae/c}
{q^{-n}e/c,\quad ae}_n,\]
which is equivalent to the binomial sum
\pnq{
&\sum_{k=0}^n
(-1)^k\binq{n}{k}q^{\binm{n-k}2}
(q^{k}a;q)_{\pav{\frac{n}{2}}}
(q^{k}c;q)_{\pav{\frac{n+1}{2}}}
\hyq{c}{q}{a,\:c}{ae,qc/e}_k\\
&\:=q^{\pav{\frac{3n^2-2n}4}}
a^{\pav{\frac{n+1}{2}}}c^{\pav{\frac{n}{2}}}
\frac{\hh{a,q/e,ae/c}{q\:}{\pav{\frac{n}{2}}}
\hh{c,e,qc/ae}{q\:}{\pav{\frac{n+1}{2}}}}
        {\hh{ae,qc/e}{q\:}{n}}.}
Now that this equation matches exactly
to \eqref{carlitz+f} specified by
\pnq{
f(k)&=q^{\pav{\frac{3k^2-2k}4}}
a^{\pav{\frac{k+1}{2}}}c^{\pav{\frac{k}{2}}}
\frac{\hh{a,q/e,ae/c}{q\:}{\pav{\frac{k}{2}}}
\hh{c,e,qc/ae}{q\:}{\pav{\frac{k+1}{2}}}}
        {\hh{ae,qc/e}{q\:}{k}},\\
g(k)&=\hyq{c}{q}{a,\:c}{ae,qc/e}_k
\quad\text{and}\quad
\vph(x;n)=(ax;q)_{\pav{\frac{n}{2}}}
(cx;q)_{\pav{\frac{n+1}{2}}};}
we have the dual relation corresponding to \eqref{carlitz+g},
which is evidenced below.
\begin{prop}[Terminating reciprocal relation]\label{pp=c}
\pnq{\hyq{c}{q}{a,\:c}{ae,qc/e}_n
=\sum_{k\ge0}\binq{n}{2k}
\frac{(1-q^{3k}c)q^{3k^2-k}(ac)^k}
     {(q^na;q)_{k}(q^nc;q)_{k+1}}
\frac{\hh{a,q/e,ae/c}{q\:}{k}
\hh{c,e,qc/ae}{q\:}{k}}
        {\hh{ae,qc/e}{q\:}{2k}}&\\
-a\sum_{k\ge0}\binq{n}{2k+1}
\frac{(1-q^{3k+1}a)q^{3k^2+2k}(ac)^k}
     {(q^na;q)_{k+1}(q^nc;q)_{k+1}}
\frac{\hh{a,q/e,ae/c}{q\:}{k}
\hh{c,e,qc/ae}{q\:}{k+1}}
        {\hh{ae,qc/e}{q\:}{2k+1}}.&}
\end{prop}
Both sums just displayed are terminating $q$-series 
with nether of them admitting closed forms. However
their combination does have a closed form.

Letting $n\to\infty$ in Proposition~\ref{pp=c} and then
applying the Weierstrass $M$-test on uniformly convergent
series,
 we get the
limiting relation:
\pnq{\hyq{c}{q}{a,\:c}{ae,qc/e}_{\infty}
=\sum_{k\ge0}
\frac{(1-q^{3k}c)q^{3k^2-k}(ac)^k}
     {(q;q)_{2k}}
\frac{\hh{a,q/e,ae/c}{q\:}{k}
\hh{c,e,qc/ae}{q\:}{k}}
        {\hh{ae,qc/e}{q\:}{2k}}&\\
-a\sum_{k\ge0}
\frac{(1-q^{3k+1}a)q^{3k^2+2k}(ac)^k}
     {(q;q)_{2k+1}}
\frac{\hh{a,q/e,ae/c}{q\:}{k}
\hh{c,e,qc/ae}{q\:}{k+1}}
        {\hh{ae,qc/e}{q\:}{2k+1}}&.}

By unifying the two sums together, we find the following theorem.
\begin{thm}[Nonterminating series identity]\label{thm=c}
\pnq{\hyq{c}{q}{a,~qc}{ae,qc/e}_{\infty}
&=\sum_{k=0}^{\infty}\frac{1-q^{3k}c}{1-c}
\frac{\hh{a,c,e,q/e,ae/c,qc/ae}{q}{k}}
     {\hh{q,ae,qc/e}{q}{2k}}
q^{3k^2-k}(ac)^k\\
&\times\bigg\{1-q^{3k}\frac{a(1-q^{3k+1}a)(1-q^kc)(1-q^ke)(1-q^{1+k}c/ae)}
	{(1-q^{3k}c)(1-q^{1+2k})(1-q^{2k}ae)(1-q^{1+2k}c/e)}\bigg\}.}
\end{thm}

Two special cases are recorded here about product of reciprocal
$q$-gamma functions, which can be utilized to establish
$q$-analogues of classical series for $\pi$ and $1/\pi$.
\begin{corl}[$a=q^{\lam}$ and $c=e=q^{1-\lam}$ in Theorem~\ref{thm=c}]\label{cc-c}
\pnq{\frac{1}{\Gam_q(\lam)\Gam_q(1-\lam)}
=&\sum_{k=0}^{\infty}q^{3 k^2}
\frac{(q^{\lam};q)^3_k(q^{1-\lam};q)^3_{k}}{(q;q)^3_{2k}}
\frac{1-q^{3k+1-\lam}}{1-q}\\
\times
&\bigg\{1-\frac{q^{3k+\lam}(1-q^{3k+1+\lam})(1-q^{k+1-\lam})^3}
      {(1-q^{3k+1-\lam})(1-q^{2k+1})^3}\bigg\}.}
\end{corl}

\begin{corl}[$a=c=q$ and $e=q^{\lam}$ in Theorem~\ref{thm=c}]\label{cc+c}
\pnq{{\Gam_q(1+\lam)\Gam_q(2-\lam)}
&=\sum_{k=0}^{\infty}\frac{1-q^{3k+1}}{1-q}
\frac{\hh{q,q,q^{\lam},q^{1-\lam},q^{\lam},q^{1-\lam}}{q}{k}}
     {\hh{q,q^{1+\lam},q^{2-\lam}}{q}{2k}}
q^{3k^2+k}\\
&\times\bigg\{1-
\frac{q^{1+3k}(1-q^{2+3k})(1-q^{1+k})(1-q^{\lam+k})(1-q^{1-\lam+k})}
     {(1-q^{1+3k})(1-q^{1+2k})(1-q^{1+\lam+2k})(1-q^{2-\lam+2k})}\bigg\}.}
\end{corl}

Five $q$-series as well as their counterparts
of classical series are displayed as follows.

\textbf{C1}. \
Recall the following series of Ramanujan~\cito{ramanujan}:
\[\frac{16}{\pi}=\sum_{k=0}^\infty
\hyp{ccc}{\frac12,\frac12,\frac12}{1,\,1,\,1}_k
\frac{5+42k}{64^k}.\]
By letting $\lam=1/2$ in Corollary~\ref{cc-c},
we recover its $q$-analogue (cf.~Chen and 
Chu~\citu{chu20RJ}{Example~40}) as follows
\pnq{\frac{1}{\Gam^2_q(\frac12)}
=\sum_{k=0}^{\infty}q^{3k^2}
\frac{(q^{1/2};q)^6_k}{(q;q)^3_{2k}}
\frac{1-q^{3k+1/2}}{1-q}
\bigg\{1-\frac{q^{3k+1/2}(1-q^{3k+3/2})}
     {(1+q^{k+1/2})^3(1-q^{3k+1/2})}\bigg\}.}

\textbf{C2}. \
For $\lam=1/4$, we get, from Corollary~\ref{cc-c},
the $q$-series identity
\pnq{\frac{1}{\Gam_q(\frac14)\Gam_q(\frac34)}
&=\sum_{k=0}^{\infty}\frac{1-q^{3k+\frac34}}{1-q}
\frac{(q^{\frac14};q)^3_k(q^{\frac34};q)^3_{k}}{(q;q)^3_{2k}}
q^{3 k^2}\\
&\times\bigg\{1-
\frac{q^{3k+\frac14}(1-q^{3k+\frac54})(1-q^{k+\frac34})^3}
      {(1-q^{3k+\frac34})(1-q^{2k+1})^3}\bigg\}.}

In order to simplify the last series, consider the series defined by
\[\sum_{k=0}^{\infty}\Lam(k),\quad\text{where}\quad
\Lam(k):=(-1)^k\frac{1-q^{\frac{1+6k}4}}{1-q}
\frac{(q^{\frac14};q^{\frac12})^3_k}{(q;q)^3_{k}}
q^{\frac34k^2}.\]
Then its \textbf{bisection series} can be reformulated as
\pnq{\sum_{k=0}^{\infty}\Lam(k)
&=\sum_{k=0}^{\infty}\big\{\Lam(2k)+\Lam(2k+1)\big\}
=\sum_{k=0}^{\infty}\Lam(2k)\Big\{1+\frac{\Lam(2k+1)}{\Lam(2k)}\Big\}\\
&=\sum_{k=0}^{\infty}\frac{1-q^{3k+\frac14}}{1-q}
\frac{(q^{\frac14};q^{\frac12})^3_{2k}}{(q;q)^3_{2k}}q^{3k^2}
\bigg\{1-
\frac{q^{3k+\frac34}(1-q^{3k+\frac74})(1-q^{k+\frac14})^3}
      {(1-q^{3k+\frac14})(1-q^{2k+1})^3}\bigg\}.}
Now it is not hard to check that
\pnq{&\frac{1-q^{3k+\frac14}}{1-q}
\bigg\{1-
\frac{q^{3k+\frac34}(1-q^{3k+\frac74})(1-q^{k+\frac14})^3}
      {(1-q^{3k+\frac14})(1-q^{2k+1})^3}\bigg\}\\
=&\frac{1-q^{3k+\frac34}}{1-q}
\bigg\{1-
\frac{q^{3k+\frac14}(1-q^{3k+\frac54})(1-q^{k+\frac34})^3}
      {(1-q^{3k+\frac34})(1-q^{2k+1})^3}\bigg\}.}
We find the following simpler series (see Chen--Chu~\citu{chu20RJ}{Example~5}
and Guo--Liu~\citu{guo18jdea}{Equation~4})
\[\frac1{\Gam_q(\frac14)\Gam_q(\frac34)}
=\sum_{k=0}^{\infty}(-1)^k\frac{1-q^{\frac{1+6k}4}}{1-q}
\frac{(q^{\frac14};q^{\frac12})^3_k}{(q;q)^3_{k}}
q^{\frac34k^2}\]
which results in a $q$-analogue of the classical one due to Guillera~\cito{gj06rj}
\[\frac{2\sqrt2}{\pi}=\sum_{k=0}^\infty\Big(\frac{-1}8\Big)^k
\hyp{ccc}{\frac12,\frac12,\frac12}{1,\,1,\,1\rule[2mm]{0mm}{2mm}}_k
\big\{1+6k\big\}.\]

\textbf{C3}. \
For $\lam=1/2$, we have, from Corollary~\ref{cc+c},
the $q$-series identity 
\pnq{\Gam^2_q\bigg(\frac32\bigg)
&=\sum_{k=0}^{\infty}q^{3 k^2+k}\frac{1-q^{3k+1}}{1-q}
\frac{(q;q)^2_k(q^{\frac12};q)^4_{k}}{(q^{\frac32};q)^2_{2k}(q;q)_{2k}}\\
&\times\bigg\{1-\frac{q^{3k+1}(1-q^{k+\frac12})(1-q^{k+1})(1-q^{3k+2})}
	{(1+q^{k+\frac12})(1-q^{2k+\frac32})^2(1-q^{3k+1})}\bigg\}}
which gives a $q$-analogue of the following series
\[\frac{9\pi}{4}=\sum_{k=0}^\infty
\hyp{cccc}{1,\frac12,\frac12,\frac12}
    {\frac54,\frac54,\frac74,\frac74\rule[2mm]{0mm}{2mm}}_k
\frac{7+42k+75k^2+42k^3}{64^k}.\]
We remark that the above $q$-series can also be derived 
by letting $x=y^2=q$ in Chu~\citu{chu18e}{Proposition~15}.

\textbf{C4}. \
Letting $a=c=e=q^{1/4}$ in Theorem~\ref{thm=c},
we derive the $q$-series identity
\pnq{\frac{\Gam_q(\frac12)}{\Gam^2_q(\frac14)}
&=\sum_{k=0}^{\infty}q^{3 k^2-\frac{k}2}\frac{1-q^{3k+\frac14}}{1-q}
\frac{(q^{\frac14};q)^4_k(q^{\frac34};q)^2_{k}}{(q^{\frac12};q)_{2k}(q;q)^2_{2k}}\\
&\times\bigg\{1-\frac{q^{3k+\frac14}(1-q^{k+\frac14})(1-q^{k+\frac34})(1-q^{3k+\frac54})}
	{(1+q^{k+\frac14})(1-q^{2k+1})^2(1-q^{3k+\frac14})}\bigg\}}
which provides a $q$-analogue of the following series
\[\frac{128\sqrt{\pi}}{\Gam^2(\frac14)}=\sum_{k=0}^\infty
\hyp{cccc}{\frac14,\frac14,\frac14,\frac34}
    {1,1,\:\frac32,\frac32\rule[2mm]{0mm}{2mm}}_k
\frac{17+396k+1392k^2+1344k^3}{64^k}.\]

\textbf{C5}. \
Letting $a=c=e=q^{3/4}$ in Theorem~\ref{thm=c},
we derive the $q$-series identity
\pnq{\frac{\Gam_q(\frac32)}{\Gam^2_q(\frac34)}
&=\sum_{k=0}^{\infty}q^{3 k^2+\frac{k}2}\frac{1-q^{3k+\frac34}}{1-q}
\frac{(q^{\frac14};q)^2_k(q^{\frac34};q)^4_{k}}{(q^{\frac32};q)_{2k}(q;q)^2_{2k}}\\
&\times\bigg\{1-\frac{q^{3k+\frac34}(1-q^{k+\frac14})(1-q^{k+\frac34})(1-q^{3k+\frac74})}
	{(1+q^{k+\frac34})(1-q^{2k+1})^2(1-q^{3k+\frac34})}\bigg\}}
which serves as a $q$-analogue of the series
\[\frac{64\sqrt{\pi}}{\Gam^2(\frac34)}=\sum_{k=0}^\infty
\hyp{cccc}{\frac14,\frac34,\frac34,\frac34}
    {1,1,\:\frac32,\frac32\rule[2mm]{0mm}{2mm}}_k
\frac{(12 k+5) (28 k+15)}{64^k}.\]

\section{Triplicate Inverse Series Relations}
For all the $n\in\mb{N}_0$, we have two equalities
\[\boxed{n=\pavv{\tfrac{1+n}3}+\pavv{\tfrac{1+2n}3}
=\pavv{\tfrac{n}3}+\pavv{\tfrac{1+n}3}+\pavv{\tfrac{2+n}3}}.\]
Then six dual relations can be established from \eqref{pfaff-q}.
However, only two of them give some interesting $q$-series identities.
Five examples are illustrated in this section without reproducing
the whole inversion procedure.

\subsection{} \
Starting from the following form of the $q$-Pfaff--Saalsch\"utz theorem \eqref{pfaff-q}
\[{_3\phi_2}\hyq{ccc}{q;q}{q^{-n},~a,~c}
{q^{-\pav{\frac{1+n}3}}ae,q^{1-\pav{\frac{2n+1}3}}c/e}
=\hyq{c}{q}{q^{-\pav{\frac{1+n}3}}e,q^{-\pav{\frac{1+n}3}}ae/c}
{q^{-\pav{\frac{1+n}3}}ae,q^{-\pav{\frac{1+n}3}}e/c}_n\]
we can derive three $q$-series identities corresponding
to the classical series of convergence rate ``$\frac4{27}$".

\textbf{D1}. \
For $a=q^{1/3}$ and $c=e=q^{2/3}$, we have the corresponding identity
\pnq{\frac{1}{\Gam_q(\frac13)\Gam_q(\frac23)}
=&\sum_{k=0}^{\infty}\frac{q^{2k^2+k}}{1-q}
\frac{\hh{q^{\frac13},q^{\frac23}}{q}{k}\hh{q^{\frac13},q^{\frac23}}{q}{2k+1}}
	{(q;q)_{k}(q;q)_{2k}(q;q)_{3k+1}}\\
\times&
\bigg\{1-
\frac{(1-q^{-k}) (1-q^{3 k+1})}
     {(1-q^{2 k+\frac{1}{3}})(1-q^{2k+\frac{2}{3}})}
+\frac{q^{2 k+1} (1-q^{k+\frac{1}{3}}) (1-q^{k+\frac{2}{3}})}
       {(1-q^{2 k+1})(1-q^{3 k+2})}\bigg\}}
which gives a $q$-analogue of the classical series
\[\frac{81\sqrt3}{2\pi}
=\sum_{k=0}^\infty\Big(\frac{4}{27}\Big)^k
\hyp{cccc}{\frac13,\frac23,\frac16,\frac56}
    {1,1,1,\:\frac32\rule[2mm]{0mm}{2mm}}_k
\big\{20+243k+414k^2\big\}.\]

\textbf{D2}. \
For $a=c=q$ and $e=q^{1/3}$, we get the corresponding identity
\pnq{\Gam_q\bigg(\frac43\bigg)\Gam_q\bigg(\frac53\bigg)
=&\sum_{k=0}^{\infty}\frac{q^{(k+1)(2k+2/3)}}{1-q}
\frac{(q^{\frac23};q)^2_{k}(q^{\frac13};q)^2_{2k+1}}
	{(q^{\frac53};q)_{k}(q^{\frac43};q)_{2k}(q;q)_{3k+1}}\\
\times&
\bigg\{1-
\frac{(1-q^{-k-\frac{2}{3}}) (1-q^{3k+1})}
     {(1-q^{2 k+\frac{1}{3}})^2}
+\frac{q^{2k+\frac43}(1-q^{k+\frac{2}{3}})^2}
      {(1-q^{2k+\frac{4}{3}})(1-q^{3 k+2})}\bigg\}}
which is a $q$-analogue of the following series
\[8\pi\sqrt3
=\sum_{k=0}^\infty\Big(\frac{4}{27}\Big)^k
\hyp{cccc}{\frac23,\frac23,\frac16,\frac16}
    {1,\frac43,\frac53,\frac76\rule[2mm]{0mm}{2mm}}_k
\big\{43+246k+414k^2\big\}.\]

\textbf{D3}. \
For $a=c=q$ and $e=q^{2/3}$, we find the corresponding identity
\pnq{\Gam_q\bigg(\frac43\bigg)\Gam_q\bigg(\frac53\bigg)
=&\sum_{k=0}^{\infty}\frac{q^{(k+1)(2k+1/3)}}{1-q}
\frac{(q^{\frac13};q)^2_{k}(q^{\frac23};q)^2_{2k+1}}
	{(q^{\frac43};q)_{k}(q^{\frac53};q)_{2k}(q;q)_{3k+1}}\\
\times&
\bigg\{1
-\frac{(1-q^{-k-\frac{1}{3}})(1-q^{3k+1})}
      {(1-q^{2k+\frac{2}{3}})^2}
+\frac{q^{2k+\frac53} (1-q^{k+\frac{1}{3}})^2}
      {(1-q^{2k+\frac{5}{3}}) (1-q^{3 k+2})}\bigg\}}
which results in a $q$-analogue of the classical series
\[40\pi\sqrt3
=\sum_{k=0}^\infty\Big(\frac{4}{27}\Big)^k
\hyp{cccc}{\frac13,\frac13,\frac56,\frac56}
    {1,\frac43,\frac53,\frac{11}6\rule[2mm]{0mm}{2mm}}_k
\big\{214+591k+414k^2\big\}.\]

\subsection{} \
Instead, rewriting the $q$-Pfaff--Saalsch\"utz theorem \eqref{pfaff-q} as
\[{_3\phi_2}\hyq{ccc}{q;q}{q^{-n},q^{\pav{\frac{n}3}}a,q^{\pav{\frac{1+n}3}}c}
{ae,q^{1-\pav{\frac{2+n}3}}c/e}
=\hyq{c}{q}{q^{-\pav{\frac{n}3}}e,q^{-\pav{\frac{1+n}3}}ae/c}
{ae,\quad q^{-\pav{\frac{2n}3}}e/c}_n\]
we obtain two further $q$-series identities.

\textbf{D4}. \
For $a=q^{1/3}$ and $c=e=q^{2/3}$, the corresponding identity reads as
\pnq{\frac{1}{\Gam_q(\frac13)\Gam_q(\frac23)}
=&\sum_{k=0}^{\infty}\frac{1-q^{4k+\frac{5}{3}}}{1-q}
\frac{(q^{\frac13};q)^2_{k}(q^{\frac23};q)^2_{k}\hh{q^{\frac13},q^{\frac23}}{q}{2k+1}}
     {(q;q)_{2k}(q;q)^2_{3k+1}}q^{5k^2+2k}\\
\times&
\bigg\{1-
\frac{(1-q^{-2 k}) (1-q^{3 k+1})^2}{(1-q^{2 k+\frac{1}{3}})
      (1-q^{2 k+\frac{2}{3}})(1-q^{4k+\frac{5}{3}})}\\
&\quad-q^{4k+\frac{4}{3}}
\frac{(1-q^{2 k+\frac{5}{3}})(1-q^{k+\frac{2}{3}})^2(1-q^{4k+\frac{7}{3}})}
     {(1-q^{2 k+1})(1-q^{3k+2})^2(1-q^{4k+\frac{5}{3}})}\bigg\}}
which provides a q-analogue of the series
\[\frac{729\sqrt3}{4\pi}
=\sum_{k=0}^\infty\Big(\frac{4}{729}\Big)^k
\hyp{cccc}{\frac13,\frac23,\frac16,\frac56}
    {1,1,1,\:\frac32\rule[2mm]{0mm}{2mm}}_k
\big\{100+1521k+2610k^2\big\}.\]

\textbf{D5}. \
For $a=c=q$ and $e=q^{1/2}$, the corresponding identity can be stated as
\pnq{\Gam^2_q\bigg(\frac32\bigg)
=&\sum_{k=0}^{\infty}\frac{q^{5k^2+\frac{3k}2}}{1+q^{\frac12}}
\frac{(q^{\frac12};q)^2_{k}(q;q)^2_{k}(q^{\frac12};q)_{2k}}
	{(q^{\frac32};q)_{3k}(q;q)_{3k}}\\
\times&\bigg\{1+q^{2 k+\frac{1}{2}}
\frac{(1-q^{2 k+\frac{1}{2}})(1-q^{4 k+2})}
     {(1-q^{3k+1})(1-q^{3 k+\frac{3}{2}})}\\
&\quad-q^{6 k+\frac{5}{2}}
\frac{(1-q^{k+\frac{1}{2}}) (1-q^{k+1})(1-q^{2k+\frac{1}{2}}) (1-q^{4 k+3})}
     {(1-q^{3 k+1}) (1-q^{3k+\frac{3}{2}})(1-q^{3 k+2}) (1-q^{3 k+\frac{5}{2}})}\bigg\}.}
By carrying on the same procedure as done for \textbf{C2},
we can confirm that the last series is, in fact,
the \textbf{bisection series} of the following one
\[\Gam^2_{q}\bigg(\frac12\bigg)
=\sum_{k=0}^{\infty}q^{\frac{k}4(3+5k)}
\frac{(q^{\frac12};q^{\frac12})^2_k(q^{\frac12};q)_k}{(q^{\frac32};q^{\frac12})_{3k}}
\frac{1+q^{\frac12+k}-q^{1+\frac{3k}2}-q^{1+2k}}{1-q^{\frac12}}.\]
This is in turn the $q$-analogue of the classical series
(cf.~Zhang~\citu{zhang2015itsf}{Example~8}): 
\[\pi=\sum_{k=0}^\infty
\Big(\frac{2}{27}\Big)^k
\hyp{cc}{1,\,\frac12}
{\frac43,\frac53\rule[2mm]{0mm}{2mm}}_k
\big(3+5k\big)
=\sum_{k=0}^\infty
\frac{6+10k}{2^k\binm{3k+2}{k+1}(k+1)(2k+1)}.\]

\textbf{Conclusive Comments.} \ We have shown that the inversion technique 
is efficient for obtaining $q$-series identities whose limiting cases result 
in $\pi$-involved series. The examples presented in this paper are far 
from exhaustive. For instance, if we start with the quadruplicate form 
of the $q$-Pfaff--Saalsch\"utz theorem \eqref{pfaff-q}
\[{_3\phi_2}\hyq{r}{q;q}{q^{-n},q^{\pav{\frac{1+n}4}}a,q^{\pav{\frac{3+n}4}}c}
{ae,q^{1-\pav{\frac{n}2}}c/e\rule[2mm]{0mm}{2mm}}
=\hyq{c}{q}{q^{-\pav{\frac{1+n}4}}e,q^{-\pav{\frac{3+n}4}}ae/c}
{ae,\qquad q^{-\pav{\frac{1+n}2}}e/c\rule[2mm]{0mm}{2mm}}_n,\]
then its dual series will give rise to the \textbf{bisection series}
of the following $q$-series
\pnq{\frac1{\Gam_q(\frac14)\Gam_q(\frac34)}
&=\sum_{k=0}^{\infty}(-1)^k
\frac{(q^{\frac14};q^{\frac12})^2_k(q^{\frac14};q^{\frac12})_{3k}}
     {(q;q)_k(q;q)^2_{2k}}q^{\frac74k^2}\\
&\times\bigg\{\frac{1-q^{\frac14+\frac{5k}2}}{1-q}
-\frac{q^{\frac34+\frac{5k}2}(1-q^{\frac14+\frac{3k}2})}
{(1-q)(1+q^{\frac14+\frac{k}2})^2(1+q^{\frac12+k})^2}\bigg\}}
which turns out to be a $q$-analogue of the elegant series
for $\frac{\sqrt2}{\pi}$ with convergence rate ``$\frac{-27}{512}$"
discovered by Guillera~\cito{gj03em}:
\[
\frac{32\sqrt2}{\pi}
=\sum_{k=0}^{\infty}\Big(\frac{-3}{8}\Big)^{3k}
\hyp{ccccc}{\frac12,\frac16,\frac56\\[-3mm]}
{1,\:1,\:1}_k\big\{15+154k\big\}.\]
We remark that the above $q$-analogue is slightly simpler
than that obtained recently by Guillera~\cito{gj18jdea}
through a totally different approach -- ``the WZ-method".


\end{document}